\input amstex\documentstyle{amsppt}
\pagewidth {12.5 cm}\pageheight {19 cm}\magnification \magstep1
\topmatter
\title Representations of reductive groups over finite rings\endtitle
\author G. Lusztig\endauthor
\address Department of Mathematics, M.I.T., Cambridge, MA 02139\endaddress
\thanks Supported in part by the National Science Foundation\endthanks
\endtopmatter   
\document  
\define\we{\wedge}
\define\cir{\circ}
\define\dw{\dot w}

\define\si{\sim}

\define\wh{\widehat}

\define\sqc{\sqcup}

\define\qua{\quad}

\define\hx{\hat x}

\define\hSi{\hat\Si}
\define\hGa{\hat\G}

\define\uT{\un T}
\define\uU{\un U}

\define\tH{\ti H}
\define\tU{\ti U}
\define\tV{\ti V}

\define\tX{\ti X}

\define\tSi{\ti\Si}

\define\op{\oplus}

\define\em{\emptyset}

\define\ra{\rangle}
\define\n{\notin}

\define\m{\mapsto}
\define\do{\dots}
\define\la{\langle}
\define\bsl{\backslash}

\define\sub{\subset}
\define\bxt{\boxtimes}
\define\T{\times}
\define\ti{\tilde}
\define\nl{\newline}
\redefine\i{^{-1}}

\define\un{\underline}

\define\ot{\otimes}
\define\bbq{\bar{\bq}_l}

\define\Ad{\text{\rm Ad}}
\define\Hom{\text{\rm Hom}}

\define\a{\alpha}
\redefine\b{\beta}
\redefine\c{\chi}
\define\g{\gamma}

\define\e{\epsilon}

\define\io{\iota}
\redefine\o{\omega}
\define\p{\pi}
\define\ph{\phi}
\define\ps{\psi}

\redefine\t{\tau}
\define\th{\theta}
\define\k{\kappa}
\redefine\l{\lambda}

\redefine\G{\Gamma}

\define\Si{\Sigma}

\define\Ph{\Phi}

\define\bF{\bold F}

\define\bn{\bold N}

\define\bq{\bold Q}

\define\bz{\bold Z}

\define\cb{\Cal B}

\define\ch{\Cal H}

\define\co{\Cal O}

\define\car{\Cal R}

\define\ct{\Cal T}

\define\cx{\Cal X}

\define\fS{\frak S}

\define\sh{\sharp}

\define\DL{DL}
\define\GE{G}
\define\LU{L}

\head Introduction\endhead
\subhead 0.1\endsubhead
In \cite{\LU, Sec.4} a cohomological construction was given (without proof)
for certain 
representations of a Chevalley group over a finite ring $R$ (arising from the ring 
of integers in a non-archimedean local field by reduction modulo a power of the 
maximal ideal); that construction was an extension of the construction of the 
virtual representations $R_T^\th$ in \cite{\DL} for groups over a finite field. One of the 
aims of this paper is to provide the missing proof. For simplicity we will assume 
that $R=\bF_{q,r}=\bF_q[[\e]]/(\e^r)$ ($\e$ is an indeterminate, $\bF_q$ is a finite
field with $q$ elements and $r\ge 1$). The general case requires only minor 
modifications. On the other hand, we treat possibly twisted groups.

Let $\bF$ be an algebraic closure of $\bF_q$. Let $G$ be a connected reductive 
algebraic group defined over $\bF$ with a given $\bF_q$-rational structure with 
associated Frobenius map $F:G@>>>G$. 

Using a cohomological method, extending that of \cite{\DL}, we will construct a 
family of irreducible representations of the finite group $G(\bF_{q,r}),r\ge 1$, 
attached to a "maximal torus" and a character of it in general position. In the case
where $r\ge 2$ and $G$ is split over $\bF_q$, the representations that we construct
are likely to be the same as those found by G\'erardin \cite{\GE} by a 
non-cohomological method (induction from a subgroup if $r$ is even; induction from a
subgroup in combination with a use of a Weil representation, if $r$ is odd, 
$\ge 3$). In any case, since for $r=1$, the cohomological construction is the only 
known construction of the generic representations, it seems natural to seek a 
cohomological construction which works uniformly for all $r\ge 1$; this is what we 
do in this paper.

In contrast with the case $r=1$, for $r\ge 2$ not all irreducible representations of
$G(\bF_{q,r})$ appear in the virtual representations that we construct. The study of
$SL_2$ with $r=2$ (see \S3) suggests that, to remedy this, one has to consider also
virtual representations attached to double cosets with respect to a "Borel subgroup"
other than those indexed by the Weyl group.

\subhead 0.2. Notation\endsubhead
Let $\e$ be an indeterminate. If $X$ is an affine algebraic variety over $\bF$ and 
$r\ge 1$, we set $X_r=X(\bF[[\e]]/(\e^r))$. Thus, if $X$ is the set of common zeroes
of the polynomials $f_i:\bF^N@>>>\bF (i=1,\do,m)$, then $X_r$ is the set of all 
$(x_1,x_2,\do,x_N)\in(\bF[[\e]]/(\e^r))^N$ such that $f_i(x_1,x_2,\do,x_N)$ (a 
priori an element of $\bF[[\e]]/(\e^r)$) is equal to $0$ for $i=1,\do,m$. We have 
$X_1=X$. For $r=0$ we set $X_r=$point. Then $X\m X_r$ is a functor from the category
of algebraic varieties over $\bF$ into itself. If $X'$ is a closed subvariety of $X$
then $X'_r$ is a closed subvariety of $X_r$. If $X$ is irreducible of dimension $d$
then $X_r$ is irreducible of dimension $dr$. For any $r\ge r'\ge 0$ we have a 
canonical morphism $\ph_{r,r'}:X_r@>>>X_{r'}$. If $r\ge 1$, we have naturally 
$X\sub X_r$ (using $\bF\sub\bF[[\e]]/(\e^r)$). If $G$ is an algebraic group over 
$\bF$ then $G_r$ is naturally an algebraic group over $\bF$. For any $r\ge r'\ge 0$,
$\ph_{r,r'}:G_r@>>>G_{r'}$ is a homomorphism of algebraic groups hence its kernel, 
$G_r^{r'}$, is a normal subgroup of $G_r$. For $r\ge 1$ we have naturally 
$G\sub G_r$. We have 
$$\{1\}=G^r_r\sub G_r^{r-1}\sub\do G_r^1\sub G_r^0=G_r.$$
For $r>r'\ge 0$, we set $G_r^{r',*}=G_r^{r'}-G_r^{r'+1}$. We have a partition
$$G_r=G_r^{0,*}\sqc G_r^{1,*}\sqc\do\sqc G_r^{r-1,*}\sqc\{1\}.$$ 
We fix a prime number $l$ invertible in $\bF$. If $X$ is an algebraic variety over
$\bF$ we write $H^j_c(X)$ istead of $H^j_c(X,\bbq)$.

For a finite group $\G$ let $\hGa=\Hom(\G,\bbq^*)$.

\subhead 0.3\endsubhead
If $T$ is a commutative algebraic group over $\bF$ with a fixed $\bF_q$-structure 
and with Frobenius map $F:T@>>>T$ we have a norm map 

$N^{F^n}_F:T^{F^n}@>>>T^F, t\m tF(t)F^2(t)\do F^{n-1}(t)$. 

\head 1. Lemmas\endhead
\proclaim{Lemma 1.1} Let $\ct,\ct'$ be two commutative, connected algebraic groups 
over $\bF$ with fixed $\bF_q$-rational structures with Frobenius maps 
$F:\ct@>>>\ct,F:\ct'@>>>\ct'$. Let $f:\ct@>\si>>\ct'$ be an isomorphism of algebraic
groups over $\bF$. Let $n\ge 1$ be such that
$F^nf=fF^n:\ct@>>>\ct'$; thus $f:\ct^{F^n}@>\si>>\ct'{}^{F^n}$. Let 
$$H=\{(t,t')\in\ct\T\ct';f(F(t)\i t)=F(t')\i t'\}.$$
(A subgroup of $\ct\T\ct'$ containing $\ct^F\T\ct'{}^F$.) Let $\th\in\wh{\ct^F}$, 
$\th'\in\wh{\ct'{}^F}$ be such that $\th\i\bxt\th'$ is trivial on 
$(\ct^F\T\ct'{}^F)\cap H^0$. Then $\th N^{F^n}_F=\th'N^{F^n}_Ff\in\wh{T^{F^n}}$. 
\endproclaim
Setting $t_1=tF(t)\do F^{n-1}(t)\in\ct,\qua t_2=f(t)F(f(t))\do F^{n-1}(f(t))\in\ct'$
for $t\in\ct$, we have
$$f(F(t_1)\i t_1)=f(tF^n(t)\i)=f(t)f(F^n(t))\i=f(t)F^n(f(t))\i=F(t_2)\i t_2,$$
so that $(t_1,t_2)\in H$. Now $t\m(t_1,t_2)$ is a morphism $\ct@>>>H$ of algebraic 
varieties and $\ct$ is connected; hence the image of this morphism is contained in 
$H^0$. In particular, if $t\in\ct^{F^n}$ we have 
$(N^{F^n}_F(t),N^{F^n}_F(f(t)))\in(\ct^F\T\ct'{}^F)\cap H^0$ hence, by assumption, 
$\th\i(N^{F^n}_F(t))\th'(N^{F^n}_F(f(t)))=1$ for all $t\in\ct^{F^n}$. The lemma is 
proved.

\subhead 1.2\endsubhead
Let $G$ be a connected reductive algebraic group over $\bF$ with a fixed
$\bF_q$-rational structure with Frobenius map $F:G@>>>G$. If $r\ge 1$ then 
$F:G@>>>G$ induces a homomorphism $F:G_r@>>>G_r$ which is the Frobenius map for a 
$\bF_q$-rational structure on $G_r$. 

Let $T,T'$ be two $F$-stable maximal tori of $G$ and let $U$ (resp. $U'$) be the 
unipotent radical of a Borel subgroup of $G$ that contains $T$ (resp. $T'$). Note 
that $U,U'$ are not necessarily defined over $\bF_q$. Let $r\ge 2$. Let 
$\ct=T_r^{r-1},\ct'=T'_r{}^{r-1}$,
$$\Si=\{(x,x',y)\in F(U_r)\T F(U'_r)\T G_r;xF(y)=yx'\}.$$
Let $N(T,T')=\{\nu\in G;\nu\i T\nu=T'\}$. Then $T$ acts on $N(T,T')$ by left
multiplication and $T'$ acts on $N(T,T')$ by right multiplication. The orbits of $T$
are the same as the orbits of $T'$; we set $W(T,T')=T\bsl N(T,T')=N(T,T')/T'$ (a 
finite set). For each $w\in W(T,T')$ we choose a representative $\dw$ in $N(T,T')$.
We have $G=\sqc_{w\in W(T,T')}G_w$ where $G_w=UT\dw U'=U\dw T'U'$. 

Let $G_{w,r}$ be the inverse image of $G_w$ under $\ph_{r,1}:G_r@>>>G$ and let
$\Si_w=\{(x,x',y)\in\Si;y\in G_{w,r}\}$. 

Now $T_r^F\T T'_r{}^F$ acts on $\Si$ by 
$(t,t'):(x,x',y)\m(txt\i,t'x't'{}\i,tyt'{}\i)$. This restricts to an action of
$T_r^F\T T'_r{}^F$ on $\Si_w$ for any $w\in W$.

If $\th\in\wh{T_r^F},\th'\in\wh{T'_r{}^F}$ and $M$ is a $T_r^F\T T'_r{}^F$-module,
we shall write $M_{\th\i,\th'}$ for the subspace of $M$ on which $T_r^F\T T'_r{}^F$
acts according to $\th\i\bxt\th'$.

\proclaim{Lemma 1.3} Assume that $r\ge 2$. Let $w\in W(T,T')$. Let 
$\th\in\wh{T_r^F}$, $\th'\in\wh{T'_r{}^F}$. Assume that 
$H^j_c(\Si_w)_{\th\i,\th'}\ne 0$ for some $j\in\bz$. Let $g=F(\dw)\i$ and let 
$n\ge 1$ be such that $g\in G^{F^n}$. Then $\Ad(g)$ carries $\ct^{F^n}$ onto
$\ct'{}^{F^n}$ and $\th|_{\ct^F}\cir N^{F^n}_F\in\wh{\ct^{F^n}}$ to
$\th'|_{\ct'{}^F}\cir N^{F^n}_F\in\wh{\ct'{}^{F^n}}$.
\endproclaim
By the definition of $G_{w,r}$, the map $U_r\T G_r^1\T(T_r\dw)\T U'_r@>>>G_{w,r}$ 
given by $(u,k,\nu,u')\m uk\nu u'$ is a locally trivial fibration with all fibres 
isomorphic to a fixed affine space. Hence the map
$$\align&
\tSi_w=\{(x,x',u,u',k,\nu)\in F(U_r)\T F(U'_r)\T U_r\T U'_r\T G_r^1\T T_r\dw;
\\&xF(u)F(k)F(\nu)F(u')=uk\nu u'x'\}@>>>\Si_w\endalign$$ 
given by $(x,x',u,u',k,\nu)\m(x,x',uk\nu u')$, is a locally trivial fibration with 
all fibres isomorphic to a fixed affine space. This map is compatible with the 
$T_r^F\T T'_r{}^F$ actions where $T_r^F\T T'_r{}^F$ acts on $\tSi_w$ by 
$$(t,t'):(x,x',u,u',k,\nu)\m(txt\i,t'x't'{}\i,tut\i,t'u't'{}\i,tkt\i,t\nu t'{}\i).
\tag a$$
Hence there exists $j'\in\bz$ such that $H^{j'}_c(\tSi_w)_{\th\i,\th'}\ne 0$. By the
substitution $xF(u)\m x,x'F(u')\i\m x'$, the variety $\tSi_w$ is rewritten as
$$\{(x,x',u,u',k,\nu)\in F(U_r)\T F(U'_r)\T U_r\T U'_r\T G_r^1\T T_r\dw;
xF(k)F(\nu)=uk\nu u'x'\};\tag b$$
in these coordinates, the action of $T_r^F\T T'_r{}^F$ is still given by (a). Let 
$$H=\{(t,t')\in\ct\T\ct';t'F(t')\i=F(\dw)\i tF(t)\i F(\dw)\}.$$
(A closed subgroup of $T_r\T T'_r$.) It acts on the variety (b) by the same formula
as in (a). (We use the fact that $hk=kh$ for any $h\in G_r^{r-1},k\in G_r^1$.) By 
\cite{\DL, 6.5}, the induced action of $H$ on $H^{j'}_c(\tSi_w)$ is trivial when 
restricted to $H^0$. In particular, the intersection $(T_r^F\T T'_r{}^F)\cap H^0$ 
acts trivially on $H^{j'}_c(\tSi_w)$. Since 
$H^{j'}_c(\tSi_w)_{\th\i,\th'}\ne 0$, it follows that $\th\i\bxt\th'$ is trivial on
$(T_r^F\T T'_r{}^F)\cap H^0$. Let $g=F(\dw)\i$ and let $n\ge 1$ be such that 
$g\in G^{F^n}$. Then $\Ad(g)$ carries $\ct^{F^n}$ onto $\ct'{}^{F^n}$ and (by Lemma
1.1 with $f=\Ad(g)$) it carries $\th|_{\ct^F}\cir N^{F^n}_F$ to 
$\th'|_{\ct'{}^F}\cir N^{F^n}_F$. The lemma is proved.

\proclaim{Lemma 1.4} Assume that $r\ge 2$. Let $\th\in\wh{T_r^F}$, 
$\th'\in\wh{T'_r{}^F}$ be such that
$$H^j_c(\Si)_{\th\i,\th'}\ne 0\tag a$$
for some $j\in\bz$. There exists $n\ge 1$ and $g\in N(T',T)^{F^n}$ such that 
$\Ad(g)$ carries $\th|_{\ct^F}\cir N^{F^n}_F\in\wh{\ct^{F^n}}$ to
$\th'|_{\ct'{}^F}\cir N^{F^n}_F\in\wh{\ct'{}^{F^n}}$.
\endproclaim
The subvarieties $G_w$ of $G$ have the following property: for some ordering $\le$ 
of $W(T,T')$, the unions $\cup_{w'\le w}G_{w'}$ are closed in $G$. It follows that 
the unions $\cup_{w'\le w}G_{w',r}$ are closed in $G_r$ and the unions 
$\cup_{w'\le w}\Si_{w'}$ are closed in $\Si$. The spectral sequence associated to 
the filtration of $\Si$ by these unions, together with (a), shows that there exists
$w\in W(T,T')$ and $j\in\bz$ such that $H^j_c(\Si_w)_{\th\i,\th'}\ne 0$. We can 
therefore apply Lemma 1.3. The lemma follows.

\subhead 1.5\endsubhead
Let $\Ph$ be the set of characters $\a:T@>>>\bF^*$ such that $\a\ne 1$ and $T$ acts
on some line $L_\a\sub Lie G$ via $\a$ (in the adjoint action); for such $\a$, let 
$G^\a$ be the one dimensional unipotent subgroup of $G$ such that $Lie G^\a=L_\a$.
For $\a\in\Ph$ there is a unique $1$-dimensional torus $T^\a$ in $T$ such that 
$T^\a$ is contained in the subgroup of $G$ generated by 
$G^\a,G^{\a\i}$. Let $\ct^\a=(T^\a)_r^{r-1}$ (a one dimensional subgroup of 
$\ct=T_r^{r-1}$).

Let $\c\in\wh{\ct^F}$. We say that $\c$ is {\it regular} if for any $\a\in\Ph$ and 
any $n\ge 1$ such that $F^n(\ct^\a)=\ct^\a$, the restriction of 
$\c\cir N^{F^n}_F:\ct^{F^n}@>>>\bbq^*$ to $(\ct^\a)^{F^n}$ is non-trivial. (It is 
enough to check that $\c\cir N^{F^n}_F|_{(\ct^\a)^{F^n}}$ is non-trivial for any 
$\a$ and for just one $n$ such that $F^n(\ct^\a)=\ct^\a$ for all $\a$.)

Let $\th\in\wh{T^F}$. We say that $\th$ is {\it regular} if $\th|_{\ct^F}$ is 
regular.

\subhead 1.6\endsubhead
Let $T$ be an $F$-stable maximal torus of $G$. Let $U,\tU,V,\tV$ be unipotent 
radicals of Borel subgroups containing $T$ such that $U\cap V=\tU\cap\tV=\{1\}$. 
Let $\Ph$ be as in 1.5. Let
$$\Ph^+=\{\a\in\Ph;G^\a\sub\tV\},\Ph^-=\{\a\in\Ph;G^\a\sub\tU\}.$$
Then $\Ph=\Ph^+\sqc\Ph^-$ and $\Ph^-=\{\a\i;\a\in\Ph^+\}$. 

For $\a\in\Ph^+$ let $ht(\a)$ be the largest integer $n\ge 1$ such that 
$\a=\a_1\a_2\do\a_n$ with $\a_i\in\Ph^+$. 

Let $x\in(G^\a)_r^b,x'\in(G^{\a'})_r^c$ where $\a,\a'\in\Ph$ and $b,c\in[0,r]$.

(a) If $b+c\ge r$ then $xx'=x'x$.

(b) If $b+c\le r$ and $\a\a'\ne 1$ then $xx'=x'xu$ where $u$ is of the form
$\prod_{i,i'\ge 1;\a^i\a'{}^{i'}\in\Ph}u_{i,i'}$ with
$u_{i,i'}\in(G^{\a^i\a'{}^{i'}})_r^{b+b'}$. 
\nl
(The factors in the last product are written in some fixed order. In the special 
case where $b+c=r-1$, these factors commute with each other by (a), since 
$r-1+r-1\ge r$.)

(c) If $b+c\ge r-1,b+2c\ge r$ and $\a\a'=1$ then $xx'=x'x\t_{x,x'}u$ where 
$\t_{x,x'}\in\ct^\a$ and $u\in(G^\a)_r^{r-1}$ are uniquely determined.

\proclaim{Lemma 1.7} We fix an order on $\Ph^+$. For any 
$z\in\tV_r,\b\in\Ph^+,$
define $x_\b^z\in G^\b_r$ by $z=\prod_{\b\in\Ph^+}x^z_\b$ (factors written using the
given order on $\Ph^+$). Let $\a\in\Ph^-,a\in[1,r-1]$. Let $z\in\tV_r^a$ be such 
that $x^z_\b\in(G^\b)_r^{a+1}$ for all $\b\in\Ph^+$ with $ht(\b)>ht(\a\i)$. Let 
$\xi\in(G^\a)_r^{r-a-1}$. Then $\xi z=z\xi\t_{\xi,z}\o_{\xi,z}$ where 
$\t_{\xi,z}\in\ct^\a$ and $\o_{\xi,z}\in\tU_r^{r-1}$.
\endproclaim
We argue by induction on $N_z=\sh(\b\in\Ph^+;x^z_\b\ne 1)$. If $N_z=0$ the result is
clear. Assume now that $N_z=1$ so that $z\in G^\b_r$ with $\b\in\Ph^+$. If $\a\b=1$,
the result follows from 1.6(c). If $\a\b\ne 1$ and $ht(\b)>ht(\a\i)$ then 
$z\in(G^\b)_r^{a+1}$ and $\xi z=z\xi$ by 1.6(b). If $\a\b\ne 1$ and 
$ht(\b)\le ht(\a\i)$ then by 1.6(b) we have $\xi z=z\xi u$ where 
$u=\prod_{i,i'\ge 1;\a^i\b^{i'}\in\Ph}u_{i,i'}$ with
$u_{i,i'}\in(G^{\a^i\b^{i'}})_r^{r-1}$; it is enough to show that if $i,i'\ge 1$, we
cannot have $\a^i\b^{i'}\in\Ph^+$. (If $\a^i\b^{i'}\in\Ph^+$ for some $i,i'\ge 1$ 
then $\a\b\in\Ph^+$ hence $ht(\b)>ht(\a\i)$, contradiction.)

Assume now that $N_z\ge 2$. We can write $z=z'z''$ where $z',z''\in\tV_r^a$, 
$N_{z'}<N_z,N_{z''}<N_z$. Using the induction hypothesis we have
$$\xi z=\xi z'z''=z'\xi\t_{\xi,z'}\o_{\xi,z'}z''$$
where $\t_{\xi,z'}\in\ct^\a$, $\o_{\xi,z'}\in\tU_r^{r-1}$. We have 
$\o_{\xi,z'}z''=z''\o_{\xi,z'}$ and $\t_{\xi,z'}z''=z''\t_{\xi,z'}$. Using again the
induction hypothesis, we have
$$\align&z'\xi\t_{\xi,z'}\o_{\xi,z'}z''=z'\xi\t_{\xi,z'}z''\o_{\xi,z'}=
z'\xi z''\t_{\xi,z'}\o_{\xi,z'}\\&=
z'z''\xi\t_{\xi,z''}\o_{\xi,z''}\t_{\xi,z'}\o_{\xi,z'}=
z\xi\t_{\xi,z'}\t_{\xi,z''}\o_{\xi,z'}\o_{\xi,z''}.\endalign$$
Thus, $\xi z=z\xi\t_{\xi,z}\o_{\xi,z}$ where 
$$\t_{\xi,z}=\t_{\xi,z'}\t_{\xi,z''},\o_{\xi,z}=\o_{\xi,z'}\o_{\xi,z''}.$$
The lemma is proved.

\subhead 1.8\endsubhead
In the setup of 1.6, let $Z=V\cap\tV$. Let $\Ph'=\{\b\in\Ph;G^\b\sub Z\}$.
We have $\Ph'\sub\Ph^+$. Let $\cx$ be the set 
of all subsets $I\sub\Ph'$ such that $I\ne\em$ and $ht:\Ph^+@>>>\bn$ is constant 
on $I$.

To any $z\in Z_r^1-\{1\}$ we associate a pair $(a,I_z)$ where $a\in[1,r-1]$ and 
$I_z\in\cx$ as follows. We define $a$ by the condition that $z\in Z_r^{a,*}$. If 
$x^z_\b\in G^\b$ are defined as in 1.8 in terms of a fixed order on $\Ph^+$, then 
$x^z_\b\in(G^\b)_r^a$ for all $\b\in\ti\Ph$ and $x^z_\b=1$ for all 
$\b\in\Ph^+-\ti\Ph$. Let $I_z$ be the set of all $\a'\in\ti\Ph$ such that 
$x_{\a'}^z\in(G^{\a'})_r^{a,*}$ and $x^z_\b\in(G^\b)_r^{a+1}$ for all $\b\in\Ph^+$
such that $ht(\b)>ht(\a')$. It is easy to see, using 1.6(a),(b), that the definition
of $I_z$ does not depend on the choice of an order on $\Ph^+$. For $a\in[1,r-1]$ and
$I\in\cx$ let $Z_r^{a,*,I}$ be the set of all $z\in Z_r^1-\{1\}$ such that 
$z\in Z_r^{a,*},I=I_z$. Thus we have a partition
$$Z_r^1-\{1\}=\sqc_{a\in[1,r-1],I\in\cx}Z_r^{a,*,I}.\tag a$$

\proclaim{Lemma 1.9}Let $T,T',U,U',r,\ct,\ct'$ be as in 1.2. Let 
$\th\in\wh{T_r^F},\th'\in\wh{T'_r{}^F}$. Assume that $\th'|_{\ct^F}=\c$ is regular. Let 
$\Si$ be as in 1.2. Then
$\sum_{j\in\bz}(-1)^j\dim H^j_c(\Si)_{\th\i,\th'}$ is equal to the number of 
$w\in W(T,T')^F$ such that $Ad(\dw):T'_r{}^F@>>>T_r^F$ carries $\th$ to $\th'$.
\endproclaim
Using the partition $\Si=\sqc_{w\in W(T,T')}\Si_w$ we see that it is enough to prove
that $\sum_{j\in\bz}(-1)^j\dim H^j_c(\Si_w)_{\th\i,\th'}$ is equal to $1$ if
$F(w)=w$ and $Ad(\dw):T'_r{}^F@>>>T_r^F$ carries $\th$ to $\th'$ and equals $0$, 
otherwise. We now fix $w\in W(T,T')$. We have 
$$\align&\Si_w\\&=\{(x,x',y)\in F(U_r)\T F(U'_r)\T G_r;xF(y)=yx',
y\in U_rG^1_r\dw T'_rU'_r=U_rZ^1_r\dw T'_rU'_r\}\endalign$$
where $Z=V\cap\dw V'\dw\i$. 
Here $V$ (resp. $V'$) is the unipotent radical of a Borel subgroup containing $T$ (resp. $T'$)
such that $U\cap V=\{1\}$ (resp. $U'\cap V'=\{1\}$. Let 
$$\align&\hSi_w
=\{(x,x',u,u',z,\t')\in F(U_r)\T F(U'_r)\T U_r\T U'_r\T Z_r^1\T T'_r;\\&
xF(u)F(z)F(\dw)F(\t')F(u')=uz\dw\t'u'x'\}.\endalign$$
The map $\hSi_w@>>>\Si_w$ given by $(x,x',u,u',z,\t')\m(x,x',uz\dw\t'u')$ is a 
locally 
trivial fibration with all fibres isomorphic to a fixed affine space. This map is 
compatible with the $T_r^F\T T'_r{}^F$-actions where $T_r^F\T T'_r{}^F$ acts on $\hSi_w$ by 
$$(t,t'):(x,x',u,u',z,\t')\m
(txt\i,t'x't'{}\i,tut\i,t'u't'{}\i,tzt\i,\dw\i t\dw\t t'{}\i).\tag a$$
Hence it is enough to show that

$\sum_{j\in\bz}(-1)^j\dim H^j_c(\hSi_w)_{\th\i,\th'}$ is equal to $1$ if $F(w)=w$ 
and $Ad(\dw):T'_r{}^F@>>>T_r^F$ carries $\th$ to $\th'$ and equals $0$, otherwise. 
\nl
By the change of variable $xF(u)\m x,x'F(u')\i\m x'$ we may rewrite $\hSi_w$ as
$$\align&\hSi_w=\{(x,x',u,u',z,\t')\in F(U_r)\T F(U'_r)\T U_r\T U'_r\T Z_r^1\T T'_r;
\\&xF(z)F(\dw)F(\t')=uz\dw\t'u'x'\}\endalign$$
with the $T_r^F\T T'_r{}^F$-action still given by (a). We have a partition 
$\hSi_w=\hSi_w'\sqc\hSi_w''$ where
$$\align&
\hSi'_w=\{(x,x',u,u',z,\t')\in F(U_r)\T F(U'_r)\T U_r\T U'_r\T(Z_r^1-\{1\})\T 
T'_r;\\&xF(z)F(\dw)F(\t')=uz\dw\t'u'x'\},\endalign$$
$$\align&\hSi''_w=\{(x,x',u,u',1,\t')\in F(U_r)\T F(U'_r)\T U_r\T U'_r\T\{1\}\T
T'_r;\\&xF(\dw)F(\t')=u\dw\t'u'x'\},\endalign$$
are stable under the $T_r^F\T T'_r{}^F$-action. It is then enough to show that

(b) $\sum_{j\in\bz}(-1)^j\dim H^j_c(\hSi''_w)_{\th\i,\th'}$ is equal to $1$ if 
$F(w)=w$ and $Ad(\dw):T'_r{}^F@>>>T_r^F$ carries $\th$ to $\th'$ and equals $0$, 
otherwise. 

(c) $H^j_c(\hSi'_w)_{\th\i,\th'}=0$ for all $j$.
\nl
We first prove (c). If $M$ is a $\ct'{}^F$-module we shall write $M_{(\c)}$ for the subspace of
$M$ on which $\ct'{}^F$ acts according to $\c$. Now $\ct'{}^F$ acts on $\hSi'_w$ by
$$t':(x,x',u,u',z,\t')\m(x,t'x't'{}\i,u,t'u't'{}\i,z,\t't'{}\i).$$
Hence $H^j_c(\hSi'_w)$ becomes a $\ct'{}^F$-module. It is enough to show that
$H^j_c(\hSi'_w)_{(\c)}=0$. 

We shall use the definitions and results in 1.6-1.8 relative to $U,\tU,V,\tV$ where
$\tU=\dw U'\dw\i$, $\tV=\dw V'\dw\i$.
 The partition 1.8(a) gives rise to a partition 
$\hSi'_w=\sqc_{a,I}\hSi^{a,I}_w$ indexed by $a\in[0,r-1],I\in\cx$ where
$$\hSi^{a,I}_w=\{(x,x',u,u',z,\t')\in\hSi'_w;z\in Z_r^{a,*,I}\}.$$
It is easy to see that there is a total order on the set of indices $(a,I)$ such 
that the union of the $\hSi^{a,I}_w$ for $(a,I)$ less than or equal than some given
$(a^0,I^0)$ is closed in $\hSi'_w$. Since the subsets $\hSi^{a,I}_w$ are stable 
under the action of $\ct'{}^F$, we see that, in order to prove (c), it is enough to 
show that
$$H^j_c(\hSi^{a,I}_w)_{(\c)}=0 \tag d$$
for any fixed $a,I$ as above. We choose $\a'\in I$. Let $\a=\a'{}\i$. Then 
$G^\a_r\sub U_r\cap \dw U'_r\dw\i$.

For any $z\in Z_r^{a,*},\xi\in(G^\a)_r^{r-a-1}$ we have 
$$\xi z=z\xi\t_{\xi,z}\o_{\xi,z}$$
where $\t_{\xi,z}\in\ct^\a,\o(\xi,z)\in\dw U'_r{}^{r-1}\dw\i$ are uniquely 
determined. (See 1.7.) Moreover, the map 
$(G^\a)_r^{r-a-1}@>>>\ct^\a,\xi\m\t(\xi,z)$ factors through
an isomorphism
$$\l_z:(G^\a)_r^{r-a-1}/(G^\a)_r^{r-a}@>\sim>>\ct^\a.$$
Let $\pi:(G^\a)_r^{r-a-1}@>>>(G^\a)_r^{r-a-1}/(G^\a)_r^{r-a}$ be the canonical 
homomorphism. We can find a morphism of algebraic varieties
$$\ps:(G^\a)_r^{r-a-1}/(G^\a)_r^{r-a}@>>>(G^\a)_r^{r-a-1}$$ 
such that $\pi\ps=1$ and $\ps(1)=1$. Let
$$\ch'=\{t'\in\ct';t'{}\i F(t')\in\dw\i\ct^\a\dw\}.$$
This is a closed subgroup of $\ct'$. For any $t'\in\ch'$ we define
$f_{t'}:\hSi^{a,I}_w@>>>\hSi^{a,I}_w$ by
$$f_{t'}(x,x',u,u',z,\t')=(xF(\xi),\hx',u,F(t')\i u'F(t'),z,\t'F(t'))$$
where 
$$\xi=\ps\l_z\i(\dw F(t')\i t'\dw\i)\in (G^\a)_r^{r-a-1}\sub U_r\cap\dw U'_r\dw\i$$
and $\hx'\in G_r$ is defined by the condition that
$$xF(\xi)F(z)F(\dw)F(\t'F(t'))=uz\dw\t'F(t')F(t')\i u'F(t')\hx'.$$
In order for this to be well defined we must check that $\hx'\in F(U'_r)$. Thus we 
must show that

$xF(\xi)F(z)F(\dw)F(\t'F(t'))\in uz\dw\t'u'F(t')F(U'_r)$
\nl
or that

$xF(z)F(\xi)F(\t_{\xi,z})F(\o_{\xi,z})F(\dw)F(\t'F(t'))\in uz\dw\t'u'F(t')F(U'_r)$.
\nl
Since $xF(z)=uz\dw\t'u'x'F(\t')\i F(\dw\i)$, it is enough to show that
$$\align&
uz\dw\t'u'x'F(\t')\i F(\dw\i)F(\xi)F(\t_{\xi,z})F(\o_{\xi,z})F(\dw)F(\t'F(t'))\\&
\in uz\dw\t'u'F(t')F(U'_r)\endalign$$ 
or that

$x'F(\t')\i F(\dw\i)F(\xi)F(\t_{\xi,z})F(\o_{\xi,z})F(\dw)F(\t'F(t'))\in 
F(t')F(U'_r)$. 
\nl
Since $x'\in F(U'_r),F(\dw\i)F(\o_{\xi,z})F(\dw)\in F(U'_r)$, it is enough to check
that

$F(\t')\i F(\dw\i)F(\xi)F(\t_{\xi,z})F(\dw)F(\t'F(t'))\in F(t')F(U'_r)$. 
\nl
Since $F(\dw\i)F(\xi)F(\dw)\in F(U'_r)$ it is enough to check that

$F(\t')\i F(\dw\i)F(\t_{\xi,z})F(\dw)F(\t'F(t'))\in F(t')F(U'_r)$
\nl
or that

$F(\dw\i)F(\t_{\xi,z})F(\dw)F(F(t'))=F(t')$
\nl
or that $\dw\i\t_{\xi,z}\dw=F(t')\i t'$ or that 
$\l_z(\p_z(\xi))=\dw F(t')\i t'\dw\i$. But this is clear.

Thus, $f_{t'}:\hSi^{a,I}_w@>>>\hSi^{a,I}_w$ is well defined for $t'\in\ch'$. It is 
clearly an isomorphism for any $t'\in\ch'$. In particular, it is a well defined 
isomorphism for any $t'\in\ch'{}^0$. By general principles, the induced map 
$f_{t'}^*:H^j_c(\hSi^{a,I}_w)@>>>H^j_c(\hSi^{a,I}_w)$ is constant when $t'$ varies 
in $\ch'{}^0$. In particular, it is constant when $t'$ varies in 
$\ct'{}^F\cap\ch'{}^0$. Now 
$\ct'{}^F\sub\ch'$ and for $t'\in\ct'{}^F$, the map $f_{t'}$ coincides with the 
action of $t'$ in the $\ct'{}^F$-action on $\hSi^{a,I}_w$. (We use that $\ps(1)=1$.)
We see that the induced action of $\ct'{}^F$ on $H^j_c(\hSi^{a,I}_w)$ is trivial 
when restricted to $\ct'{}^F\cap\ch'{}^0$. 

We can find $n\ge 1$ such that $F^n(\dw\i\ct^\a\dw)=\dw\i\ct^\a\dw$. Then
$t'\m t'F(t')F^2(t')\do F^{n-1}(t')$ is a well defined morphism 
$\dw\i\ct^\a\dw@>>>\ch'$. Its image is a connected subgroup of $\ch'$ hence is 
contained in $\ch'{}^0$. If $t'\in(\dw\i\ct^\a\dw)^{F^n}$ then 
$N^{F^n}_F(t')\in\ct'{}^F$; thus, $N^{F^n}_F(t')\in\ct'{}^F\cap\ch'{}^0$. We see 
that the action of $N^{F^n}_F(t')\in\ct'{}^F$ on $H^j_c(\hSi^{a,I}_w)$ is trivial 
for any $t'\in(\dw\i\ct^\a\dw)^{F^n}$. 

If we assume that $H^j_c(\hSi^{a,I}_w)_{(\c)}\ne 0$, it follows that 
$t'\m\c(N^{F^n}_F(t'))$ is the trivial character of $(\dw\i\ct^\a\dw)^{F^n}$. This 
contradicts our assumption that $\c$ is regular. Thus, (d) holds. Hence (c) holds.

We now prove (b). Let
$$\tH=\{(t,t')\in T_r\T T'_r;tF(t)\i=F(\dw)t'F(t')\i F(\dw\i)\}.$$
This is a closed subgroup of $T_r\T T'_r$ containing $T_r^F\T T'_r{}^F$. Now the action of 
$T_r^F\T T'_r{}^F$ on $\hSi_w''$ extends to an action of $\tH$ given by the same 
formula. To see this consider $(t,t')\in\tH$ and $(x,x',u,u',1,\t')\in\hSi''_w$. We
must show that 

$(txt\i,t'x't'{}\i,tut\i,t'u't'{}\i,1,\dw\i t\dw \t't'{}\i)\in\hSi''_w$
\nl
that is,

$txt\i F(\dw)F(\dw\i)F(t)F(\dw)F(\t')F(t'{}\i)=tut\i\dw\dw\i t\dw\t't'{}\i 
t'u't'{}\i t'x't'{}\i$
\nl
or that

$xt\i F(t)F(\dw)F(\t')F(t'{}\i)=u\dw\t'u'x't'{}\i$
\nl
or that
$xt\i F(t)F(\dw)F(\t')F(t'{}\i)=xF(\dw)F(\t')t'{}\i$
\nl
or that $t\i F(t)F(\dw)F(t'{}\i)=F(\dw)t'{}\i$;
this is clear. Let $T_*,T'_*$ be the reductive part of $T_r,T'_r$ (thus $T_*$ is a 
torus isomorphic to $T$). Let $\tH_*=\tH\cap(T_*\T T'_*)$. Then $\tH_*^0$ is a torus
acting on $\hSi_w''$ by restriction of the $\tH$-action. The fixed point set 
$(\hSi_w'')^{\tH_*^0}$ of the $\tH_*^0$-action is stable under the action of 
$T_r^F\T T'_r{}^F$ and by general principles we have
$$\sum_{j\in\bz}(-1)^j\dim H^j_c(\hSi''_w)_{\th\i,\th'}=
\sum_{j\in\bz}(-1)^j\dim H^j_c((\hSi''_w)^{\tH^0_*})_{\th\i,\th'}.$$
It is then enough to show that

(e) $\sum_{j\in\bz}(-1)^j\dim H^j_c((\hSi''_w)^{\tH^0_*})_{\th\i,\th'}$ is equal to
$1$ if $F(w)=w$ and $Ad(\dw):T'_r{}^F@>>>T_r^F$ carries $\th$ to $\th'$ and equals 
$0$, otherwise. 
\nl
Let $(x,x',u,u',1,\t)\in(\hSi''_w)^{\tH_*^0}$. By Lang's theorem the first 
projection $\tH_*@>>>T_*$ is surjective. It follows that the first projection 
$\tH_*^0@>>>T_*$ is surjective. Similarly the second projection $\tH_*^0@>>>T'_*$ is
surjective. Hence for any $t\in T_*,t'\in T'_*$ we have 

$txt\i=x, t'x't'{}\i=x',tut\i=u,t'u't'{}\i=u'$
\nl
hence $x=x'=u=u'=1$. Thus, $(\hSi''_w)^{\tH_*^0}$ is contained in 

(f) $\{(1,1,1,1,1,\t');\t'\in T'_r,F(\dw\t')=\dw\t'\}$. 
\nl
The set (f) is clearly contained in the fixed point set of $\tH$. Note that (f) is 
empty unless $F(w)=w$. We can therefore assume that $F(w)=w$. In this case, (f) is 
stable under the action of $\tH$. In particular it is stable under the action of
$\tH_*^0$. Since (f) is finite and $\tH_*^0$ is connected, we see that $\tH_*^0$ 
must act trivially on (f). Thus, (f) is exactly the fixed point set of $\tH_*^0$. 
Hence this fixed point can be identified with 
$(\dw T'_r)^F$. From this (e) follows easily. The lemma is proved.

\head 2. The main results\endhead
\subhead 2.1\endsubhead
Let $G,F$ be as in 1.2. Let $T$ be an $F$-stable maximal torus in $G$ and let $U$ be
the unipotent radical of a Borel subgroup of $G$ that contains $T$. (Note that $U$ 
is not necessarily $F$-stable.) Let $r\ge 1$. Let $\car(G_r^F)$ be the group of
virtual representations of $G_r^F$ over $\bbq$. Let $\la,\ra$ be the standard inner
product $\car(G_r^F)\T\car(G_r^F)@>>>\bz$. Let 
$$S_{T,U}=\{g\in G_r;g\i F(g)\in F(U_r)\}.$$
The finite group $G_r^F\T T_r^F$ acts on $S_{T,U}$ by $(g_1,t):g\m g_1gt\i$. For any
$i\in\bz$ we have an induced action of $G_r^F\T T_r^F$ on $H^i_c(S_{T,U})$. For
$\th\in\wh{T^F_r}$, we denote by $H^i_c(S_{T,U})_\th$ the subspace of 
$H^i_c(S_{T,U})$ on which $T_r^F$ acts according to $\th$. This is a 
$G^F_r$-submodule of $H^i_c(S_{T,U})$. Let
$$R_{T_r,U_r}^\th=\sum_{i\in\bz}(-1)^iH^i_c(S_{T,U})_\th\in\car(G_r^F).$$

\proclaim{Proposition 2.2} Assume that $r\ge 2$. Let $(T',U',\th')$ be another 
triple like $T,U,\th$. Let $\ct=T_r^{r-1},\ct'=T'_r{}^{r-1}$. 

(a) Let $i,i'$ be integers. Assume that there exists an irreducible $G_r^F$-module 
that appears in the $G_r^F$-module $(H^i_c(S_{T,U})_{\th\i})^*$ (dual of 
$H^i_c(S_{T,U})_{\th\i}$) and in the $G_r^F$-module $H^{i'}_c(S_{T',U'})_{\th'}$.
There exists $n\ge 1$ and $g\in N(T',T)^{F^n}$ such that $\Ad(g)$ carries
$\th\cir N^{F^n}_F|_{\ct^{F^n}}\in\wh{\ct^{F^n}}$ to 
$\th'\cir N^{F^n}_F|_{\ct'{}^{F^n}}\in\wh{\ct'{}^{F^n}}$.

(b) Assume that there exists an irreducible $G_r^F$-module that appears in the 
virtual $G_r^F$-module $\sum_i(-1)^iH^i_c(S_{T,U})_\th$ and in the virtual
$G_r^F$-module $\sum_i(-1)^iH^i_c(S_{T',U'})_{\th'}$. There exists $n\ge 1$ and 
$g\in N(T',T)^{F^n}$ such that $\Ad(g)$ carries 
$\th\cir N^{F^n}_F|_{\ct^{F^n}}\in\wh{\ct^{F^n}}$ to 
$\th'\cir N^{F^n}_F|_{\ct'{}^{F^n}}\in\wh{\ct'{}^{F^n}}$.
\endproclaim  
We prove (a). Consider the free $G_r^F-$action on $S_{T,U}\T S_{T',U'}$ given by 
$g_1:(g,g')\m(g_1g,g_1g')$. The map 
$$(g,g')\m(x,x',y),x=g\i F(g),x'=g'{}\i F(g'),y=g\i g'$$
defines an isomorphism of $G_r^F\bsl(S_{T,U}\T S_{T',U'})$ onto $\Si$ (as in 1.2).

The action of $T^F_r\T T'_r{}^F$ on $S_{T,U}\T S_{T',U'}$ given by right 
multiplication by $t\i$ on the first factor and by $t'{}\i$ on the second factor 
becomes an action of $T_r^F\T T'_r{}^F$ on $\Si$ given by 
$(x,x',y)\m(txt\i,t'x't'{}\i,tyt'{}\i)$. Our assumption implies that the 
$G_r^F$-module $H^i_c(S_{T,U})_{\th\i}\ot H^{i'}_c(S_{T',U'})_{\th'}$ contains the 
unit representation with non-zero multiplicity. Hence the subspace of 
$H^{i+i'}_c(G_r^F\bsl(S_{T,U}\T S_{T',U'}))$ on which $T_r^F\T T'_r{}^F$ acts 
according to $\th\i\bxt\th'$ is non-zero. Equivalently, 
$H^{i+i'}_c(\Si)_{\th\i,\th'}\ne 0$. We now use Lemma 1.4; (a) follows.

We prove (b). By general principles we have
$$\sum_i(-1)^i(H^i_c(S_{T,U})_{\th\i})^*=\sum_i(-1)^iH^i_c(S_{T,U})_\th.$$ 
Hence the assumption of (b) implies that the assumption of (a) holds. Hence the 
conclusion of (a) holds. The proposition is proved.

\proclaim{Proposition 2.3}We preserve the setup of 2.2. Assume that $\th$ or $\th'$
is regular. The inner product $\la R_{T_r,U_r}^\th,R_{T'_r,U'_r}^{\th'}\ra$ is 
equal to the number of $w\in W(T,T')^F$ such that $Ad(\dw):T'_r{}^F@>>>T_r^F$ 
carries $\th$ to $\th'$.
\endproclaim
We may assume that $\th'$ is regular. As in the proof of 2.2, we have
$$\align\la R_{T_r,U_r}^\th,R_{T'_r,U'_r}^{\th'}\ra\\&=\sum_{i,i'\in\bz}
(-1)^{i+i'}\dim(H^i_c(S_{T,U})_{\th\i}\ot H^{i'}_c(S_{T',U'})_\th)^{G_r^F}\\&
=\sum_{j\in\bz}(-1)^j\dim H^j_c(G_r^F\bsl(S_{T,U}\T S_{T',U'}))_{\th\i,\th'}\\&
=\sum_{j\in\bz}(-1)^j\dim H^j_c(\Si)_{\th\i,\th'}\endalign$$
where $()^{G_r^F}$ denotes the space of $G_r^F$-invariants.
It remains to use 1.9.

\proclaim{Corollary 2.4} Assume that $r\ge 2$. Let $T,U$ be as in 2.1. Assume that 
$\th\in\wh{T_r^F}$ is regular. 

(a) $R_{T_r,U_r}^\th$ is independent of the choice of $U$.

(b) Assume also that the stabilizer of $\th$ in $W(T,T)^F$ is $\{1\}$. Then 
$R_{T_r,U_r}^\th$ is $\pm$ an irreducible $G_r^F$-module.
\endproclaim
We prove (a). Let $U'$ be the unipotent radical of another Borel subgroup of $G$ 
containing $T$. Let $R=R_{T_r,U_r}^\th$, $R'=R_{T_r,U'_r}^\th$. By 2.3 we have 

$\la R,R\ra=\la R,R'\ra=\la R',R\ra=\la R',R'\ra$. 
\nl
Hence $\la R-R',R-R'\ra=0$, so that $R=R'$. This proves (a). 
In the setup of (b), we see from 2.3 that
$\la R_{T_r,U_r}^\th,R_{T_r,U_r}^\th\ra=1$. This proves (b).

\subhead 2.5\endsubhead
Assume that $r\ge 2$. Let $T$ be as in 2.1. Assume that $\th\in\wh{T_r^F}$ is 
regular. We set 

$R_{T_r}^\th=R_{T_r,U_r}^\th$
\nl
where $U$ is chosen as in 2.1. (By 2.4(a), this is independent of the choice of 
$U$.)

\head 3. An example\endhead
\subhead 3.1\endsubhead
Let $A=\bF[[\e]]/(\e^2)$. Define $F:A@>>>A$ by $F(a_0+\e a_1)=a_0^q+\e a_1^q$ where
$a_0,a_1\in\bF$. Let $V$ be a $2$-dimensional $\bF$-vector space with a fixed 
$\bF_q$-rational structure with Frobenius map $F:V@>>>V$. Let $G=SL(V)$. Then $G$ 
has an $\bF_q$-rational structure with Frobenius map $F:G@>>>G$ such that 
$F(gv)=F(g)F(v)$ for all $g\in G,v\in V$. Let $V_2=A\ot_{\bF}V$. Then $G_2$ (see 
0.2) may be identified with the group of all automorphisms of the free $A$-module 
$V_2$ with determinant $1$. We regard $V$ as a subset of $V_2$ by $v\m 1\ot v$. Any
element of $V_2$ can be written uniquely in the form $v_0+\e v_1$ where 
$v_0,v_1\in V$. The Frobenius map $F:V_2@>>>V_2$ satisfies 
$F(v_0+\e v_1)=F(v_0)+\e F(v_1)$ for $v_0,v_1\in V$.

Let $\wh{G_2^F}$ be the set of isomorphism classes of irreducible representations of
$G_2^F$ over $\bbq$. The objects of $\wh{G_2^F}$ can be classified using the fact 
that $G_2^F$ is a semidirect product of $G^F$ and $\bF_q^3$.

The table below shows the number of representations in $\wh{G_2^F}$ of various 
dimensions assuming that $q$ is odd; the first column indicates the dimension, the 
second column indicates the number of representations of that dimension. 
$$\matrix  \dim &\sh\\
            1  & 1 \\
           q&   1 \\
          q+1 &   (q-3)/2\\
          (q+1)/2&    2  \\
           q-1 &    (q-1)/2\\
          (q-1)/2&    2\\                  
           q^2+q &   (q-1)^2/2\\
           q^2-q & (q^2-1)/2   \\
           (q^2-1)/2& 2q \\               \endmatrix$$
The analogous table in the case where $q$ is a power of $2$ is
$$\matrix  \dim &\sh\\
            1  & 1 \\
           q&   1 \\
          q+1 &   (q-2)/2\\
          q-1 &    q/2\\
           q^2+q &   (q-1)(q-2)/2\\
           (q^2+q)/2 & 2(q-1)\\
           q^2-q & (q^2-q)/2   \\
           (q^2-q)/2 & 2(q-1)   \\
           q^2-1& q \\               \endmatrix$$
\subhead 3.2\endsubhead
Let $\cb$ be the set of all $A$-submodules $L\sub V_2$ such that $L$ is a direct
summand of $V_2$ and $L$ is free of rank $1$. Now $G_2$ acts transitively on $\cb$.
The diagonal action on $\cb\T\cb$ has three orbits $\co,\co',\co''$ where

$\co=\{(L,L')\in\cb\T\cb;L=L'\}$,

$\co'=\{(L,L')\in\cb\T\cb;L\cap L'=\e L=\e L'\}$,

$\co''=\{(L,L')\in\cb\T\cb;L\cap L'=0\}$.
\nl
If $L\in\cb$ then $F(L)\in\cb$.
Thus we obtain a map $F:\cb@>>>\cb$, the Frobenius map of a  
$\bF_q$-rational structure on $\cb$. Let 

$X=\{L\in\cb;(L,F(L))\in\co\}$, 

$X'=\{L\in\cb;(L,F(L))\in\co'\}$, 

$X''=\{L\in\cb;(L,F(L))\in\co''\}$.
\nl
Then $X,X',X''$ form a partition of $\cb$ into $G_2^F$-stable subvarieties.

We now define some finite coverings of $X,X',X''$ as follows. Let $e,e'$ be an 
$\bF$-basis of $V$ such that $F(e)=e,F(e')=e'$. Let $\uT$ be the subgroup of $G$ 
consisting of the automorphisms $e\m ae,e'\m a\i e'$ with $a\in\bF^*$. (An $F$-stable maximal torus of 
$G$.) Let $\uU$ be the subgroup of $G$ consisting of 
the automorphisms 

$e\m e+be',e'\m e'$ with $b\in\bF^*$. 
\nl
Let $\nu\in G$ be such that $\nu(e)=e',\nu(e')=-e$. Let $h\in G$ be such that $he=e,he'=e'+\e e$. Let 

$\tX=\{g\in G_2;g\i F(g)\in\uU_2\}/\uU_2$, 

$\tX'=\{g\in G_2;g\i F(g)\in h\uU_2\}/(\uU_2\cap h\uU_2h\i)$,

$\tX''=\{g\in G_2;g\i F(g)\in\nu\uU_2\}$.
\nl
(We use the action of $\uU_2$ or $\uU_2\cap h\uU_2h\i$ on $G_2$ by right 
translation.) Then $g\m Age'$ is a well defined morphism 
$\tX@>>>X,\tX'@>>>X',\tX''@>>>X''$. This is a finite principal covering with group 
$\G,\G',\G''$ respectively (acting by right translation) where

$\G=\uT_2^F$ (of order $q^2-q$) 

$\G'=\{x\in\uT_2\uU_2; x\i hF(x)\in h\uU\}/(\uU\cap h\uU h\i)\cong\{\pm 1\}\T\bF_q$
(of order $2q$ if $q$ is odd, of order $q$ if $q$ is a power of $2$) 

$\G''=\{t\in\uT_2;F(t)=t\i\}$ (of order $q^2+q$).

For any variety $Y$ with an action of a finite group and any character $\o$ of that
finite group, let $H^j_c(Y)_\o$ denote the subspace of $H^j_c(Y)$ on which the 
finite group acts according to $\o$. Thus, for $\o$ in $\hGa$ (resp.
$\hGa',\hGa''$), $H^j_c(\tX)_\o$ (resp. $H^j_c(\tX')_\o,H^j_c(\tX'')_\o$) is well 
defined.

\subhead 3.3\endsubhead
Let $Y=\{g\in G;g\i F(g)\in h\uU_2\}$. We wish to describe $Y$ more explicitly. If 
$g\in G_2$, the condition that $g\in Y$ is that 

$F(g)e=ghue=gh(e+xe')=g(e+xe'+\e xe)$, 

$F(g)e'=ghue'=gh(e')=g(e'+\e e)$
\nl
for some $x\in A$. Define $a,b,c,d\in A$ by $ge=ae+be',ge'=ce+de'$. The condition 
that $g\in Y$ is

$F(a)e+F(b)e'=ae+be'+\e xae+\e xbe'+xce+xde'$, 

$F(c)e+F(d)e'=ce+de'+\e ae+\e be'$
\nl
for some $x\in A$. Thus, we may identify $Y$ with the set of all $(a,b,c,d)\in A^4$
such that

$F(a)=a+\e xa+xc,F(b)=b+\e xb+xd,F(c)=c+\e a,F(d)=d+\e b,ad-bc=1$
\nl
for some $x\in A$, or equivalently, such that

$(F(a)-a)(\e b+d)=(F(b)-b)(\e a+c),F(c)=c+\e a,F(d)=d+\e b,ad-bc=1$.
\nl
Setting $a=a_0+\e a_1,b=b_0+\e b_1,c=c_0+\e c_1,d=d_0+\e d_1$ with $a_i,d_i\in\bF$,
we see that $Y$ is identified with the set consisting of all
$(a_0,b_0,c_0,d_0,a_1,b_1,c_1,d_1)\in\bF^8$ such that

(a) $c^q_0=c_0,d^q_0=d_0,c^q_1=c_1+a_0,d^q_1=d_1+b_0$,

(b) $a_0d_0-b_0c_0=1,a_0d_1+a_1d_0-b_0c_1-b_1c_0=0$,

(c) $(a^q_0-a_0)d_0=(b^q_0-b_0)c_0$, 

$(a^q_0-a_0)b_0+(a^q_0-a_0)d_1+(a^q_1-a_1)d_0
=(b^q_0-b_0)a_0+(b^q_0-b_0)c_1+(b^q_1-b_1)c_0$.
\nl
Actually the equations (c) are a consequence of the other equations, hence they can
be omitted. The first equation (b) can be written (using (a)):

$(c^q_1-c_1)d_0-(d^q_1-d_1)c_0=1$,
\nl
that is,

$(c_1d_0-d_1c_0)^q-(c_1d_0-d_1c_0)=1$.
\nl
Setting $f=c_1d_0-d_1c_0$, we see that $Y$ is identified with the set of all

$(a_0,b_0,c_0,d_0,a_1,b_1,c_1,d_1,f)\in\bF^9$
\nl
such that

$c^q_0=c_0,d^q_0=d_0,c^q_1=c_1+a_0,d^q_1=d_1+b_0$,

$f^q-f=1, c_1d_0-d_1c_0=f, a_0d_1+a_1d_0-b_0c_1-b_1c_0=0$.
\nl
Now on $Y$ we have a free right action of $\uU_2^1=\uU_2\cap h\uU_2h\i$, $u:g\m gu$.
In terms of coordinates, this is $(a,b,c,d)\m(a+\e xc,b+\e xd,c,d), x\in\bF$ or
$$(a_0,b_0,c_0,d_0,a_1,b_1,c_1,d_1,f)\m
(a_0,b_0,c_0,d_0,a_1+xc_0,b_1+xd_0,c_1,d_1,f).$$
The set of orbits $Y/\uU_2^1=\tX'$ may be identified with the set of all

$(a_0,b_0,c_0,d_0,c_1,d_1,f)\in\bF^7$
\nl
such that
$$c^q_0=c_0,d^q_0=d_0,c^q_1=c_1+a_0,d^q_1=d_1+b_0,f^q-f=1,c_1d_0-d_1c_0=f.$$
We consider the obvious projection of this set to the finite set 

$\{(c_0,d_0,f)\in\bF^3;c^q_0=c_0,d^q_0=d_0,f^q-f=1,(c_0,d_0)\ne(0,0)\}$.
\nl
The fibre of this projection at $(c_0,d_0,f)\in\bF^3$ is the affine line 
$\{(c_1,d_1)\in\bF^2;c_1d_0-d_1c_0=f\}$. Thus, $\tX'$ is a union of $(q^2-1)q$ 
affine lines. Hence $H^j_c(\tX')=0$ for $j\ne 2$ and $H^2_c(\tX')$ is a permutation
representation of $G_2^F$ of dimension $(q^2-1)q$. For $q$ odd, it follows easily 
that, as a $G_2^F$-module, $H^2_c(\tX')$ is the direct sum of the $2q$ irreducible 
representations of degree $(q^2-1)/2$ (each one with multiplicity $1$); more 
precisely, for any $\o'\in\hGa'$, $H^2_c(\tX')_{\o'}$ is irreducible of degree 
$(q^2-1)/2$ and each irreducible representation  of degree $(q^2-1)/2$ is obtained 
for exactly one $\o'$.

Similarly, for $q$ a power of $2$, $H^2_c(\tX')$ is the direct sum of the $q$ irreducible 
representations of degree $q^2-1$ (each one with multiplicity $1$); more 
precisely, for any $\o'\in\hGa'$, $H^2_c(\tX')_{\o'}$ is irreducible of degree 
$q^2-1$ and each irreducible representation  of degree $q^2-1$ is obtained for 
exactly one $\o'$.

\subhead 3.4\endsubhead
Now $\tX$ is a permutation representation of $G_2^F$ of dimension $q^4-q^2$ which is
easy to analyze. We see that $H^j_c(\tX)=0$ for $j\ne 0$ and, for $q$ odd,
$H^0_c(\tX)$ is the direct sum of all irreducible representations of degree $q^2+q$
and $q+1$ (each one with multiplicity $2$), those of degree $1,q,(q+1)/2$ (each one
with multiplicity $1$) and $H^2_c(\tX')_{\o'}$ with $\o'\in\hGa'$, $\o'{}^2=1$ (each
one with multiplicity $q-1$). More precisely, if $\o\in\hGa$ then $H^0_c(\tX)_\o$ is

irreducible of degree $q^2+q$ if $\o|_{\uT_2^1{}^F}\ne 1$;

the direct sum of $\op_{\o'\in\hGa;\o'{}^2=1}H^2_c(\tX')_{\o'}$ with an 
irreducible representations of degree $q+1$, if $\o|_{\uT_2^1{}^F}=1,\o^2\ne 1$;

the direct sum of $\op_{\o'\in\hGa;\o'{}^2=1}H^2_c(\tX')_{\o'}$ with the two
irreducible representations of degree $(q+1)/2$, if 
$\o|_{\uT_2^1{}^F}=1,\o^2=1,\o\ne 1$;

the direct sum of $\op_{\o'\in\hGa;\o'{}^2=1}H^2_c(\tX')_{\o'}$ with the two
irreducible representations of degree $1$ and $q$, if $\o=1$.

On the other hand, for $q$ a power of $2$, $H^0_c(\tX)$ is the direct sum of all irreducible 
representations of degree $q^2+q$ and $q+1$ (each one with multiplicity $2$), those
of degree $1,q,(q^2+q)/2$ (each one with multiplicity $1$) and $H^2_c(\tX')_1$ 
(with multiplicity $q-1$). More precisely, if $\o\in\hGa$ then $H^0_c(\tX)_\o$ is 

irreducible of degree $q^2+q$ if $\o|_{\uT_2^1{}^F}\ne 1,\o^2\ne 1$;

the direct sum of two irreducible representations of degree $(q^2+q)/2$ if 
$\o^2=1,\o\ne 1$;

the direct sum of $H^2_c(\tX')_1$ with an irreducible representations of degree 
$q+1$, if $\o|_{\uT_2^1{}^F}=1,\o\ne 1$;

the direct sum of $H^2_c(\tX')_1$ with the two irreducible representations of degree
$1$ and $q$, if $\o=1$.

\subhead 3.5\endsubhead
Let 
$$\fS_0=\{x_0\in V;x_0\we F(x_0)=e\we e'\},\quad \fS_{00}=\{x_0\in\fS_0;F^2(x_0)=-x_0\}.$$
Now $G^F$ acts on $\fS_0$ (restriction of the $G$-action on $V$). This restricts to
a $G^F$-action on $\fS_{00}$. We show that this action is simply transitive. If 
$g\in G^F$ keeps fixed some $x_0\in\fS_{00}$ then it also keeps fixed $F(x_0)$ hence
it must be $1$ (recall that $x_0,F(x_0)$ form a basis of $V$). Thus the $G^F$-action
on $\fS_{00}$ has trivial isotropy. We may identify $\fS_{00}$ with 
$\{(a,b)\in\bF^2;ab^q-a^qb=1,a^{q^2}=-a,b^{q^2}=-b\}$. For such $(a,b)$ we have 
automatically $a\ne 0$. We make a change of variable $(a,b)\m(a,c)$ where $c=b/a$. 
Then $\fS_{00}$ becomes 

$\{(a,c)\in\bF^2;a^{q+1}(c^q-c)=1,a^{q^2}=-a,c^{q^2}=c\}$.
\nl
The second projection maps this to $\{c\in\bF;c^{q^2}=c,c^q\ne c\}$ which has 
$q^2-q$ elements. The fibre at $c$ is $\{a\in\bF;a^{q+1}=(c^q-c)\i\}$. (For such $a$
we have automatically $a^{q^2}=-a$ since $c^{q^2}=c$.) This fibre has exactly $q+1$
elements since $(c^q-c)\i\ne 0$. We see that $\sh(\fS_{00})=(q+1)(q^2-q)=\sh(G^F)$. 
It follows that the $G^F$-action on $\fS_{00}$ is indeed simply transitive.

\subhead 3.6\endsubhead
We now analyze $\tX''$. Let
$$\fS=\{x\in V_2;x\we F(x)=e\we e'\}.$$
Now $G_2^F$ acts on $\fS$ by $g_1:x\m g_1x$. The map $g\m g(e')$ defines an 
isomorphism 

$\io:\tX''@>\si>>\fS$.
\nl
We check that this is a well defined bijection. Let $g\in\tX''$. Then $F(g)=g\nu u$ 
for some $u\in\uU_2$. Let $x=ge'$. Then for some $u\in\uU_2$ we have
$$\align&x\we F(x)=(ge')\we F(ge')=(ge')\we F(g)e'=e'\we g\i F(g)e'=e'\we\nu ue'\\&
=e'\we\nu e'=e'\we(-e)=e\we e',\endalign$$
hence $x\in\fS$ and $\io$ is well defined. Now let $x\in\fS$. We can find $g\in G_2$
such that $ge'=x$. Then 
$$e\we e'=x\we F(x)=(ge')\we F(ge')=(ge')\we F(g)e'=e'\we g\i F(g)e'.$$
Hence $g\i F(g)e'=-e+be'$ for some $b\in A$. It follows that $g\i F(g)=u'\nu u$ 
where $u,u'\in\uU_2$. Then $(gu')\i F(gu')=\nu uF(u')$ hence $gu'\in\tX''$. Clearly,
$\io(gu')=x$ so that $\io$ is surjective. Now assume that $g,g'\in\tX''$ satisfy 
$\io(g)=\io(g')$ that is $ge'=g'e'$. Then $g'=gu'$, $u'\in\uU_2$. We have 
$g'{}\i F(g')=\nu u$ with $u\in\uU_2$ hence $u'{}\i g\i F(g)F(u')=\nu u$. Also, 
$g\i F(g)=\nu\ti u$ with $\ti u\in\uU_2$ hence $u'{}\i\nu\ti uF(u')=\nu u$ so that
$u'\in\nu\uU_2\nu\i$. Thus, $u'\in\uU_2\cap(\nu\uU_2\nu\i)=\{1\}$ hence $u'=1$ and 
$g'=g$. Thus, $\io$ is injective hence bijective. It commutes with the 
$G_2^F$-actions.

Now $\fS$ consists of the elements $x_0+\e x_1$, with $x_0,x_1\in V$ such that

$(x_0+\e x_1)\we(F(x_0)+\e F(x_1))=e\we e'$,
\nl
that is 

$x_0\we F(x_0)=e\we e'$ and $x_1\we F(x_0)+x_0\we F(x_1)=0$. 
\nl
We have a morphism 

$\k:\fS@>>>\fS_0,x_0+\e x_1\m x_0$. 
\nl
If $x_0\in\fS_0$ then $\k\i(x_0)$ may be identified with

$\{x_1\in V;x_1\we F(x_0)+x_0\we F(x_1)=0\}$.
\nl
Note that $x_0,F(x_0)$ form a basis of $V$ hence $F^2(x_0)=c_0x_0+c_1F(x_0)$ with
$c_0,c_1\in\bF$. Since $x_0\we F(x_0)=e\we e'$ is $F$-stable, we have
$x_0\we F(x_0)=F(x_0)\we F^2(x_0)$ hence $c_0=-1$. Let $\fS_{01}=\fS_0-\fS_{00}$.
We have a partition $\fS=\fS_*\cup\fS_{**}$ where 
$\fS_*=\k\i(\fS_{00}),\fS_{**}=\k\i(\fS_{01})$ are $G_2^F$-stable. If $x_0\in\fS_0$,
then any $x_1\in V$ can be written uniquely in the form 

$x_1=a_0x_0+a_1F(x_0)$
\nl
with $a_0,a_1\in\bF$. The condition that $x_0+\e x_1\in\k\i(x_0)$ is

$(a_0x_0+a_1F(x_0))\we F(x_0)+x_0\we(a_0^qF(x_0)+a_1^qF^2(x_0))=0$,
\nl
that is,

$a_0x_0\we F(x_0)+x_0\we(a_0^qF(x_0)-a_1^qx_0+a_1^qc_1F(x_0))=0$,
\nl
that is,

$a_0x_0\we F(x_0)+a_0^qx_0\we F(x_0)+a_1^qc_1x_0\we F(x_0)=0$,
\nl

or
$a_0+a_0^q+a_1^qc_1=0$.
\nl
Thus we may identify $\k\i(x_0)$ with $\{(a_0,a_1)\in\bF^2;a_0+a_0^q+a_1^qc_1=0\}$.
If $c_1\ne 0$ (that is if $x_0\in\fS_{01}$) this is isomorphic to the affine line. 
Thus, $\k$ restricts to an affine line bundle $\fS_{**}@>>>\fS_{01}$.

Now the action of $\G''$ on $\tX''$ corresponds under $\io$ to the action of 
$\{\l\in A;\l F(\l)=1\}$ on $\fS$ by scalar multiplication. Hence the action of 
$\{t\in\G'';t\in T_2^1\}$ on $\tX''$ corresponds to the action of 
$A'=\{\l\in A;\l F(\l)=1,\l\in 1+\e A\}$ on $\fS$ by scalar multiplication. The 
action of $1+\e\l_1\in A'$ (with $\l_1\in\bF$) in the coordinates $(x_0,a_0,a_1)$ is
$(x_0,a_0,a_1)\m(x_0,a_0+\l_1,a_1)$. Thus it preserves each fibre of $\k$.

Now $\fS_{**}$ is stable under the action of $\{\l\in A;\l F(\l)=1\}$ and the 
restriction of this action to $A'$ preserves each fibre of $\fS_{**}@>>>\fS_{01}$ 
(an affine line); hence this group acts trivially on $H^j_c()$ of each such fibre 
hence it also acts trivially on $H^j_c(\fS_{**})$. Thus, 
$H^j_c(\fS)@>>>H^j_c(\fS_*)$ is an isomorphism on the part where $\sum_{\l\in A'}\l$
acts as $0$.

We now study $H^j_c(\fS_*)$. If $x_0\in\fS_{00}$ then $\k\i(x_0)$ may be identified 
with $\{(a_0,a_1)\in\bF^2;a_0+a_0^q=0\}$. Thus, $\fS_*$ is an affine line bundle 
over 

$\fS_{00}\T\{a_0\in\bF;a_0+a_0^q=0\}$
\nl
which is a transitive permutation 
representation of $G_2^F$ that is explicitly known from 3.5. It follows that 
$H^j_c(\fS_*)=0$ for $j\ne 2$ and the part of $H^2_c(\fS_*)$ where 
$\sum_{\l\in A'}\l$ acts as $0$ is the direct sum of the irreducible representations
of degree $q^2-q$ (each one with multiplicity $2$) and of degree $(q^2-q)/2$ (each 
one with multiplicity $1$); note that the latter representations occur only when $q$
is a power of $2$.

We now study the part of $H^j_c(\fS)$ where $A'$ acts as $1$. This is the same as 
$H^j_c(A'\bsl\fS)$. The map $(x_0,a_0,a_1)\m(x_0,\ti a_0,a_1),\ti a_0=a_0+a_0^q$ is
an isomorphism of $A'\bsl\fS$ with the set of all 
$(x_0,\ti a_0,a_1)\in\fS_0\T\bF\T\bF$ such that $\ti a_0+a_1^qc_1=0$. (Here $c_1$ is
determined by $x_0$ as above.) Hence the map $(x_0,a_0,a_1)\m(x_0,a_1)$ is an 
isomorphism $A'\bsl\fS@>\si>>\fS_0\T\bF$. Thus, $H^j_c(A'\bsl\fS)=H^{j-2}_c(\fS_0)$. 
Thus, $G_2^F$ acts on $H^j_c(A'\bsl\fS)$ through its quotient $G^F$ and that action
is explicitly known from the representation theory of $G^F$.

We see that $H^4_c(\tX'')$ is the $1$ dimensional representation; $H^3_c(\tX'')$ is
the direct sum of all irreducible representations of degree $q-1$ (each one with 
multiplicity $2$) and those of degree $(q-1)/2,q$ (each one with multiplicity $1$);
$H^2_c(\tX'')$ is the direct sum of all irreducible representations of degree 
$q^2-q$ (each one with multiplicity $2$) and of degree $(q^2-q)/2$ (each one with 
multiplicity $1$); $H^j_c(\tX'')=0$ for $j\n\{2,3,4\}$; note that the 
representations of degree $(q-1)/2$ occur only for $q$ odd, while those of degree
$(q^2-q)/2$ occur only for $q$ a power of $2$.

More precisely, if $\o\in\hGa''$ and $q$ is odd, then 

$H^4_c(\tX'')_\o$ is irreducible of degree $1$ if $\o=1$ and is $0$ otherwise;

$H^3_c(\tX'')_\o$ is irreducible of degree $q-1$ if $\o|_{\G''\cap T_2^1}=1$,
$\o^2\ne 1$; it is the direct sum of two irreducible representations of degree
$(q-1)/2$ if $\o|_{\G''\cap T_2^1}=1$, $\o^2=1,\o\ne 1$; it it is irreducible of 
degree $q$ if $\o=1$; it is $0$ if $\o|_{\G''\cap T_2^1}\ne 1$;

$H^2_c(\tX'')_\o$ is irreducible of degree $q^2-q$ if $\o|_{\G''\cap T_2^1}\ne 1$
and is $0$ otherwise.

Similarly, if $\o\in\hGa''$ and $q$ is a power of $2$, then 

$H^4_c(\tX'')_\o$ is irreducible of degree $1$ if $\o=1$ and is $0$ otherwise;

$H^3_c(\tX'')_\o$ is irreducible of degree $q-1$ if $\o|_{\G''\cap T_2^1}=1$,
$\o\ne 1$; it is irreducible of degree $q$ if $\o=1$; it is $0$ if 
$\o|_{\G''\cap T_2^1}\ne 1$;

$H^2_c(\tX'')_\o$ is irreducible of degree $q^2-q$ if $\o|_{\G''\cap T_2^1}\ne 1$,
$\o^2\ne 1$; it is the direct sum of two irreducible representations of degree
$(q^2-q)/2$ if $\o^2=1,\o\ne 1$; it is $0$ otherwise.

\subhead 3.7\endsubhead
We see that any irreducible representation of $G_2^F$ appears in at least one of the
representations $H^j_c(\tX)_\o,H^j_c(\tX')_\o,H^j_c(\tX'')_\o$. More precisely, the
regular representation of $G_2^F$ is a $\bq$-linear combination of the virtual
representations
$$\sum_{j\in\bz}(-1)^jH^j_c(\tX)_\o,\sum_{j\in\bz}(-1)^jH^j_c(\tX')_\o,
\sum_{j\in\bz}(-1)^jH^j_c(\tX'')_\o.$$ 

\subhead 3.8\endsubhead
Let $\g\in G$ be such that $\g\i F(\g)=\nu$. We set $T=\g\uT\g\i$, $U=\g\uU\g\i$. 
Then $T$ is an $F$-stable maximal torus of $G$ and $U$ is the unipotent radical of a
Borel subgroup of $G$ containing $T$. Hence $S_{T,U}$ is defined (with $r=2$). Now
$g\m g\g\i$ defines an isomorphism 

$\tX''@>\si>>S_{T,U}$
\nl
and an isomorphism $\G''@>\si>>T^F_2$. Also $G_2^F\T\G''$ acts on $\tX''$ by 
$(g_1,t):g\m g_1gt\i$. This action is compatible with the $G_2^F\T T_2^F$-action on
$S_{T,U}$ via the isomorphisms above. We see that the virtual representations 
$\sum_{j\in\bz}(-1)^jH^j_c(\tX'')_\o$ of $G_2^F$ are the same as the virtual 
representations $R_{T,U}^\th$.

\subhead 3.9\endsubhead
We return to the general setup of 1.2. Let $\un B$ be an $F$-stable Borel subgroup
of $G$ with unipotent radical $\uU$. For any $x\in G_r$ let 
$X_x=\{g\in G_r;g\i F(g)\in x\uU_r\}$. Then $G_r^F$ acts on $X_x$ by left
translations hence it acts naturally on $H^j_c(X_x)$. Note that the isomorphism 
class of the $G_r^F$-module $H^j_c(X_x)$ depends only on the $(\uU_r,\uU_r)$ double
coset of $x$ in $G_r$. We conjecture that any irreducible representation of $G_r^F$
appears in $\sum_{j\in\bz}(-1)^jH^j_c(X_x)$ for some $x\in G_r$. This holds for 
$r=1$ (see \cite{\DL}) and also for $G=SL_2,r=2$, by the results in this section.

\enddocument